# Partial Semigroup Algebras Associated to Partial Actions

B. Tabatabaie Shourijeh[1] and S. Moayeri Rahni [2]

**Abstract**

For a given inverse semigroup $S$, we introduce the notion of algebraic crossed product by using a given partial action of S, and we will prove that under some condition it is associative. Also we will introduce the concept of partial semigroup algebra $K_{Par}(S)$, and we show that the suitable quotient of $K_{Par}(S)$ is a sort of crossed product.

**Keywords:** Inverse semigroup, Partial action, Partial Representation, Crossed Product.

## Introduction

Let $A$ and $B$ be Banach algebras and $\alpha$ be a multiplicative linear functional The notion of a *partial crossed product* of a $C^*$-algebra by the group of integers was defined by R. Exel in [2]. Roughly, the automorphism used in the definition of a crossed product of a $C^*$-algebra $A$ by the group of integers was replaced by an isomorphism between two ideals of $A$, namely, partial automorphisms. K. McCalanahan in [5] defined a partial crossed product of a $C^*$-algebra by a discrete group. The ideas involved were mostly straightforward generalization of the definitions given in [2]. K. McCalanahan used partial action of groups to define a partial crossed products. This motivates us to define the algebraic crossed products by a partial actions of an inverse semigroups $S$, and we will discuss the associativity of this construction. To define a partial semigroup algebra, we need some definitions and terminologies that A. Buss and R. Exel defined in [1]. The two important concepts are a *universal inverse semigroup* $Pr(S)$ and *partial action* of the inverse semigroups. We will use the universal inverse semigroup $Pr(S)$ to define a partial semigroup algebra $Pr(S)$. Here, we will assume that the reader is familiar to the inverse semigroups and the concepts of partial actions of inverse semigroups, multiplier algebras, and partial representation. For a thorough treatment of these subjects the reader is referred to [1, 3, 4, 5, 7]. In section 1, we will introduce the *algebraic crossed product*, and we will discuss the conditions under which this construction is associative. In section 2, to prove our main result, the suitable quotient of the partial semigroup algebra is a sort of crossed products, we need to define some partial action. To do this, we will show that we can define a partial representation for an inverse semigroup $S$ from a partial action of it and vise versa. At the end, we will prove that the suitable quotient of $K_{Par}(S)$ is a sort of crossed product.

1) Department of Mathematics, College of Sciences, Shiraz University, Shiraz, Iran.
   E-mail: tabataba@math.susc.ac.ir
2) Department of Mathematics, College of Sciences, Shiraz University, Shiraz, Iran.
   E-mail: smoayeri@shirazu.ac.ir

# Results and Discussion

## 1. Algebraic Crossed Products

Throughout this work, we will assume that $S$ is a unital inverse semigroup, $\mathsf{A}$ is an associative algebra, and $\alpha$ is a partial action of $S$ on $\mathsf{A}$.

Let $\alpha = (\{\alpha_s\}_{s \in S}, \{X_s\}_{s \in S})$ be a partial action of $S$ on $\mathsf{A}$. Let $L = \{\sum_{s \in S} a_s \delta_s : a_s \in X_s\}$, the set of all formal finite sums, with the following multiplication:
$$(a_s \delta_s).(b_t \delta_t) = \alpha_s(\alpha_{s^*}(a_s) b_t) \delta_{st}.$$

Note that $\alpha_s(\alpha_{s^*}(a_s) b_t)$ is an element of $X_{st}$ simply because
$$\alpha_s(\alpha_{s^*}(a_s) b_t) \in \alpha_s(X_{s^*} \cap X_t) = X_s \cap X_{st}.$$

Thus the multiplication is well-defined, and $L$ is an algebra with this multiplication. Let $I$ be the ideal of $L$ generated by the set $\{a\delta_r - a\delta_t : \text{where } r \leq t \text{ and } a \in X_r\}$. Note that by [1, Proposition 3.8] $X_r \subseteq X_t$ since $r \leq t$. Now, we are at a position to introduce the *algebraic crossed product* $\mathsf{A} \times_\alpha S$ associated to a partial action $\alpha$ as we proposed in the beginning. the *algebraic crossed product* $\mathsf{A} \times_\alpha S$ is the algebra $\dfrac{L}{I}$. Considering the quotient map $\Phi : L \to \mathsf{A} \times_\alpha S$ from $L$ onto $\mathsf{A} \times_\alpha S$, we will denote the element $\Phi(a_s \delta_s) \in \mathsf{A} \times_\alpha S$ by $\overline{a_s \delta_s}$.

With the aid of the following pivotal Theorem, we are able to answer the associativity question. We will show that the mentioned multiplication is associative if for each $s \in S$ the ideals $X_s$ associated to partial action $\alpha$ are $(L,R)$-associative, that is, when $L \circ R' = R' \circ L$ for all multipliers $(L, R), (L', R')$ of $X_s$.

**Theorem 1.1** If $\alpha$ is a partial action of a unital inverse semigroup $S$ on an algebra $\mathsf{A}$ such that each $X_s$ $(s \in S)$ is $(L,R)$-associative, then the algebraic crossed product $\mathsf{A} \times_\alpha S$ is associative.

*Proof.* The Proof is a slight modification of [3, Theorem 3.4].
As a consequence of Theorem 1.1, we have the following corollary.

**Corollary 1.2** Let $\alpha$ be a partial action of a unital inverse semigroup $S$ on an algebra $\mathsf{A}$ such that each $X_s$ is an idempotent or non-degenerate ideal *of* $\mathsf{A}$ for each $s \in S$. Then $\mathsf{A} \times_\alpha S$ is associative.
*Proof.* This follows from the fact that $(L, R)$-associativity is equivalent to the condition of the Corollary and Theorem 1.1.

**Definition 1.3** We shall say that an algebra $\mathsf{A}$ is strongly associative if for any unital inverse semigroup $S$ and any partial action of $S$ on $\mathsf{A}$, namely $\alpha$, the algebraic crossed product $\mathsf{A} \times_\alpha S$ is associative.

As a consequence of Corollary 1.2 we have:

**Corollary 1.4** A semiprime algebra $\mathsf{A}$ is strongly associative.

## 2. Partial Representation

Throughout this section, by $\pi$ we mean a partial representation of an inverse semigroup $S$ on a unital algebra $\mathsf{B}$. For simplicity hereafter, we will denote an element $\overline{a_s \delta_s}$ of $\mathsf{A} \times_\alpha S$ by $a_s \delta_s$.

**Lemma 2.1** Suppose that $\alpha$ is a partial action of $S$ on an algebra $\mathsf{A}$, where $X_s$ (the rang of $\alpha_s$) is a unital algebra with the unit element $1_s$. Then the map $\pi_\alpha : S \ni s \mapsto 1_s \delta_s \in \mathsf{A} \times_\alpha S$ is a partial representation.

*Proof.* First of all, note that since each $X_s$ is unital, it is idempotent, hence, each $X_s$ is $(L,R)$-associative, thus, $\mathsf{A} \times_\alpha S$ is associative by Theorem 1.1. To prove the Lemma, let $s,t$ are elements of $S$. Then

$$\pi_\alpha(s^*)\pi_\alpha(s)\pi_\alpha(t) = 1_{s^*}\delta_{s^*}.1_s\delta_s.1_t\delta_t$$
$$= \alpha_{s^*}(1_s.1_s)\delta_{s^*s}.1_t\delta_t$$
$$= (1_{s^*}1_{s^*s})\delta_{s^*s}.1_t\delta_t. \tag{1}$$

The last equality is obtained from the fact that $1_s 1_t$ is the unit element of the algebra $X_s \cap X_t$. Since $\alpha$ is a partial action, the equality $\alpha_s(X_{s^*} \cap X_t) = X_s \cap X_{st}$ implies that

$$\alpha_s(1_{s^*}1_t) = 1_s 1_{st}. \tag{2}$$

Now, we back to the proof. By the definition of the multiplication, the right hand side of (1) is equal to

$$\alpha_{s^*}(\alpha_{s^*}(1_{s^*}1_{s^*s})1_t)\delta_{s^*st} = \alpha_{s^*}(1_{s^*}1_{s^*s}1_t)\delta_{s^*st}. \tag{3}$$

Note that $X_{s^*} \subseteq X_{s^*s}$ because

$$X_{s^*} = \alpha_{s^*}(X_s \cap X_s) = X_{s^*} \cap X_{s^*s} \subseteq X_{s^*s}.$$

Thus, $1_{s^*}1_{s^*s} = 1_{s^*s}1_{s^*} = 1_{s^*}$. The foregoing computations and equation (2) gives that the right hand side of (3) is equal to

$$1_{s^*}1_{s^*st}\delta_{s^*st}. \tag{4}$$

On the other hand,

$$\pi_\alpha(s^*)\pi_\alpha(st) = 1_{s^*}\delta_{s^*}.1_{st}\delta_{st}$$
$$= \alpha_{s^*}(\alpha_s(1_{s^*})1_{st})\delta_{s^*st}$$
$$= \alpha_{s^*}(1_s 1_{st})\delta_{s^*st}$$
$$= 1_{s^*}1_{s^*st}\delta_{s^*st}. \tag{5}$$

Comparing (4) and (5), one has that $\pi_\alpha(s^*)\pi_\alpha(s)\pi_\alpha(t) = \pi_\alpha(s^*)\pi_\alpha(st)$. For $s, t \in S$, we have

$$\pi_\alpha(s)\pi_\alpha(t)\pi_\alpha(t^*) = 1_s \delta_s . 1_t \delta_{tt^*} = \alpha_s(1_{s^*} 1_t)\delta_{stt^*}$$

$$= 1_s 1_{st} 1_{stt^*} \delta_{stt^*}, \quad (6)$$

by using (2) and the fact that $1_s 1_{st} \in X_{stt^*}$. On the one hand, we have that

$$\pi_\alpha(st)\pi_\alpha(t^*) = 1_{st} \delta_{st} . 1_{t^*} \delta_{t^*}$$
$$= \alpha_{st}(1_{(st)^*} 1_{t^*})\delta_{stt^*}$$
$$= 1_{st} 1_{stt^*} \delta_{stt^*}.$$

The last equality is obtained by using (2). Since $stt^* \leq s$, we have that $X_{stt^*} \subseteq X_s$ by [1, Proposition 3.8], therefore, $1_{st} 1_{stt^*} \delta_{stt^*} = 1_s 1_{st} 1_{stt^*} \delta_{stt^*}$. Thus, the proof is complete if we compare this equality with (6).

**Definition 2.2** Tow partial representations $\pi : S \mapsto \mathsf{B}$ and $\pi' : S \mapsto \mathsf{B}'$ are equivalent if there exists an isomorphism of algebras $\phi : \mathsf{B} \mapsto \mathsf{B}'$ such that
$$\pi(s) = \phi(\pi'(s))$$
for all $s \in S$.

**Remark 1** It is easy to see that if $\alpha$ and $\alpha'$ are two equivalent partial actions for $S$ then $\pi_\alpha$ and $\pi_{\alpha'}$ are equivalent partial representations of $S$.

Our next goal is to define a partial action from a given partial representation of an inverse semigroup $S$. To do this, let $\pi : S \mapsto \mathsf{B}$ be a partial representation of $S$ into a $\mathsf{K}$-algebra $\mathsf{B}$. Considering the elements $\varepsilon_s = \pi(s)\pi(s^*)$ for $s \in S$, we have the following Lemma.

**Lemma 2.3** Let $s, t$ are in $S$. Then
  (i) For given $s, t$ in $S$, $\varepsilon_s$ and $\varepsilon_t$ are commuting and idempotent elements,
  (ii) and $\pi(s)\varepsilon_t = \varepsilon_{st}\pi(s)$.

*Proof.* Taking $s \in S$, one has that
$$\varepsilon_s \varepsilon_s = \pi(s)\pi(s^*)\pi(s)\pi(s^*) = \pi(s)\pi(s^*) = \varepsilon_s,$$
thus $\varepsilon_s$ is idempotent. On the one hand, for $s, t \in S$ we have

$$\varepsilon_s \varepsilon_t = \pi(s)\pi(s^*)\pi(t)\pi(t^*)$$
$$= \pi(s)\pi(s^* t)\pi(t^*)$$
$$= \pi(s)\pi(s^* t)\pi(t^* s)\pi(s^* t)\pi(t^*)$$
$$= \pi(ss^* t)\pi(t^* s)\pi(s^* tt^*). \quad (7)$$

Interchanging the role of $s$ and $t$, from the foregoing argument we have that

$$\varepsilon_t \varepsilon_s = \pi(tt^*s)\pi(s^*t)\pi(t^*ss^*)$$
$$= \pi(tt^*s)\pi(s^*t)\pi(t^*s)\pi(s^*t)\pi(t^*ss^*)$$
$$= \pi(tt^*ss^*t)\pi(t^*s)\pi(s^*tt^*ss^*)$$
$$= \pi(ss^*t)\pi(t^*s)\pi(s^*tt^*) = \varepsilon_s \varepsilon_t.$$

For part (ii), let $s, t \in S$, then

$$\pi(s)\varepsilon_t = \pi(s)\pi(t)\pi(t^*)$$
$$= \pi(s)\pi(s^*)\pi(s)\pi(t)\pi(t^*)$$
$$= \pi(s)\pi(t)\pi(t^*)\pi(s)\pi(s^*)$$
$$= \pi(st)\pi(t^*s^*)\pi(s) = \varepsilon_{st}\pi(s). \tag{8}$$

Let $\mathsf{A}$ be the subalgebra of $\mathsf{B}$ generated by all $\varepsilon_s$ ($s \in S$), and for a fixed $s \in S$ set $X_s = \varepsilon_s \mathsf{A}$. Now, we can provide a partial action from $\pi$ as follows:

**Lemma 2.4** The maps $\alpha_s^\pi : X_{s^*} \mapsto X_s$ ($s \in S$) defined by $\alpha_s^\pi(a) = \pi(s)a\pi(s^*)$ ($a \in X_{s^*}$) are isomorphisms of $\mathsf{K}$-algebras which determine a partial action $\alpha^\pi$ of $S$ on the algebra $\mathsf{A}$.

***Proof.*** For simplicity, let $\alpha^\pi = \alpha$. We show that $\mathsf{A}$ is invariant under the map $a \mapsto \pi(s)a\pi(s^*)$ for each $s \in S$. Since $\mathsf{A}$ is spanned by elements of the form $\varepsilon_{t_1}...\varepsilon_{t_n}$ with $t_1...t_n \in S$. For such elements, one has that

$$\pi(s)\varepsilon_{t_1}...\varepsilon_{t_n}\pi(s^*) = \varepsilon_{st_1}...\varepsilon_{st_n}\pi(s)\pi(s^*)$$
$$= \varepsilon_{st_1}...\varepsilon_{st_n}\varepsilon_s \in \mathsf{A}$$

Thus, $\pi(s)\mathsf{A}\pi(s^*) \subseteq \mathsf{A}$. Now, for $s \in S$,

$$\alpha_s(\varepsilon_{s^*}) = \pi(s)\varepsilon_{s^*}\pi(s^*)$$
$$= \pi(s)\pi(s^*)\pi(s)\pi(s^*)$$
$$= \varepsilon_s \varepsilon_s = \varepsilon_s.$$

Consequently, for the elements of the form $\varepsilon_{s^*}\varepsilon_{t_1}...\varepsilon_{t_n}$, one has that

$$\pi(s)\varepsilon_{s^*}\varepsilon_{t_1}...\varepsilon_{t_n}\pi(s^*) = \varepsilon_{ss^*}\varepsilon_{st_1}...\varepsilon_{st_n}\pi(s)\pi(s^*)$$
$$= \varepsilon_s \varepsilon_{ss^*}\varepsilon_{st_1}...\varepsilon_{st_n} \in X_s.$$

Since $X_{s^*}$ is spanned by the element of the form $\varepsilon_{s^*}\varepsilon_{t_1}...\varepsilon_{t_n}$ with $t_1...t_n \in S$, we have that $\alpha_s : X_{s^*} \mapsto X_s$. Moreover, it is a homomorphism of algebras. Indeed, taking elements $a = \varepsilon_{s^*}\varepsilon_{t_1}...\varepsilon_{t_n}$ and $b = \varepsilon_{s^*}\varepsilon_{r_1}...\varepsilon_{r_n}$ in $X_{s^*}$, we see that

$$\alpha_s(a)\alpha_s(b) = \pi(s)\varepsilon_{s^*}\varepsilon_{t_1}...\varepsilon_{t_n}\pi(s^*)\pi(s)\varepsilon_{s^*}\varepsilon_{r_1}...\varepsilon_{r_n}\pi(s^*)$$

$$= \pi(s)\varepsilon_{s^*}\varepsilon_{t_1}...\varepsilon_{t_n}\varepsilon_s\varepsilon_{s^*}\varepsilon_{r_1}...\varepsilon_{r_n}\pi(s^*)$$
$$= \pi(s)ab\pi(s^*) = \alpha_s(ab).$$

Obviously, for $s \in S$, $\varepsilon_s$ is the identity of $X_s$ and $\alpha_s(\varepsilon_{s^*}) = \varepsilon_s$, thus, $\alpha_s$ is a homomorphism of algebras. Observe that for $s \in S$ we have $\alpha_s(\alpha_{s^*}(\varepsilon_s)) = \varepsilon_s$. By induction, we shall prove that $\alpha_s \circ \alpha_{s^*}(\varepsilon_s\varepsilon_{t_1}...\varepsilon_{t_n}) = \varepsilon_s\varepsilon_{t_1}...\varepsilon_{t_n}$. For $n = 1$,

$$\alpha_s \circ \alpha_{s^*}(\varepsilon_s\varepsilon_{t_1}) = \alpha_s(\pi(s^*)\varepsilon_s\varepsilon_{t_1}\pi(s))$$
$$= \alpha_s(\varepsilon_{s^*s}\varepsilon_{s^*s}\varepsilon_{s^*t_1}\varepsilon_{s^*s})$$
$$= \pi(s)\varepsilon_{s^*s}\varepsilon_{s^*s}\varepsilon_{s^*t_1}\pi(s^*)$$
$$= \pi(s)\pi(s^*)\pi(s)\pi(s^*s)\pi(s^*s)\pi(s^*t_1)\pi(t_1^*s)\pi(s^*)$$
$$= \pi(s)\pi(s^*)\pi(ss^*s)\pi(s^*s)\pi(s^*t_1)\pi(t_1^*s)\pi(s^*)$$
$$= \pi(s)\pi(s^*)\pi(s)\pi(s^*)\pi(t_1)\pi(t_1^*)\pi(s)\pi(s^*)$$
$$= \varepsilon_s\varepsilon_{t_1}.$$

Now, suppose that for $n < k$ the statement is true. Let $\varepsilon_s\varepsilon_{t_1}...\varepsilon_{t_k}$ be any element of $X_s$. Then

$$\alpha_s \circ \alpha_{s^*}(\varepsilon_s\varepsilon_{t_1}...\varepsilon_{t_k}) = \alpha_s \circ \alpha_{s^*}(\varepsilon_s\varepsilon_{t_1}...\varepsilon_{t_{k-1}}\varepsilon_s\varepsilon_{t_k})$$
$$= \alpha_s \circ \alpha_{s^*}(\varepsilon_s\varepsilon_{t_1}...\varepsilon_{t_{k-1}})\alpha_s \circ \alpha_{s^*}(\varepsilon_s\varepsilon_{t_k}),$$

so, the induction hypothesis implies that

$$\alpha_s \circ \alpha_{s^*}(\varepsilon_s\varepsilon_{t_1}...\varepsilon_{t_k}) = \varepsilon_s\varepsilon_{t_1}...\varepsilon_{t_k}.$$

Since $X_s$ is spanned by such elements, we deduce that $\alpha_s \circ \alpha_{s^*}$ is the identity map on $X_s$. Interchanging the role of $s$ and $s^*$, we see that $\alpha_{s^*} \circ \alpha_s$ is the identity map on $X_{s^*}$, hence, $\alpha_s^{-1} = \alpha_{s^*}$, thus, $\alpha_s$ is an isomorphism. Let $s,t \in S$ and write an element $a \in X_{s^*} \cap X_t$ as $a = \varepsilon_{s^*}\varepsilon_t b$ with $b \in A$. So,

$$\alpha_s(a) = \pi(s)\varepsilon_{s^*}\varepsilon_t b\pi(s^*)$$
$$= \varepsilon_{ss^*}\varepsilon_{st}\pi(s)b\pi(s^*)$$
$$= \varepsilon_{st}\varepsilon_{ss^*}\varepsilon_{st}\pi(s)b\pi(s^*) \in X_{st},$$

where in the last step we have used the fact that $\pi(s)A\pi(s^*) \subseteq A$. This gives $\alpha_s(X_{s^*} \cap X_t) \subseteq X_s \cap X_{st}$. Taking $x \in X_s \cap X_{st}$, write $x$ as $x = \varepsilon_s\varepsilon_{st}b$ with $b \in A$. Moreover, we can assume that $b = \varepsilon_{t_1}...\varepsilon_{t_n}$ with $t_1...t_n \in S$. Note that

$$\varepsilon_{st}\varepsilon_s = \pi(st)\pi(t^*s^*)\pi(s)\pi(s^*)$$
$$= \pi(st)\pi(t^*)\pi(s^*)\pi(s)\pi(s^*)$$
$$= \pi(s)\pi(t)\pi(t^*)\pi(s^*)\pi(s)\pi(s^*)$$
$$= \alpha_s(\varepsilon_t\varepsilon_{s^*}).$$

On the other hand, we have that

$$x = \varepsilon_s \varepsilon_{st} \varepsilon_{t_1} \ldots \varepsilon_{t_n}$$
$$= \varepsilon_s \varepsilon_{st} \pi(s) \pi(s^*) \varepsilon_s \varepsilon_{t_1} \ldots \varepsilon_{t_n}$$
$$= \varepsilon_s \varepsilon_{st} \pi(s) \varepsilon_{s^*s} \varepsilon_{s^*s} \varepsilon_{s^*t_1} \ldots \varepsilon_{s^*t_n} \pi(s^*)$$
$$= \alpha_s(\varepsilon_t \varepsilon_{s^*}) \alpha_s(\varepsilon_{s^*s} \varepsilon_{s^*s} \varepsilon_{s^*t_1} \ldots \varepsilon_{s^*t_n})$$
$$= \alpha_s(\varepsilon_t \varepsilon_{s^*} \varepsilon_{s^*s} \varepsilon_{s^*s} \varepsilon_{s^*t_1} \ldots \varepsilon_{s^*t_n}).$$

From the fact that $\varepsilon_t \varepsilon_{s^*} \varepsilon_{s^*s} \varepsilon_{s^*s} \varepsilon_{s^*t_1} \ldots \varepsilon_{s^*t_n} \in X_{s^*} \cap X_t$, we have $x \in \alpha_s(X_{s^*} \cap X_t)$, and we deduce that $\alpha_s(X_{s^*} \cap X_t) = X_s \cap X_{st}$. Let $a \in \alpha_t^{-1}(X_t \cap X_{s^*}) \subseteq X_{t^*} \cap X_{t^*s^*}$. Then $\varepsilon_{t^*} a = a \varepsilon_{t^*} = a$, so,

$$\alpha_s \circ \alpha_t(a) = \alpha_s(\pi(t) a \pi(t^*)) = \pi(s) \pi(t) a \pi(t^*) \pi(s^*)$$
$$= \pi(s) \pi(t) a \varepsilon_{t^*} \pi(t^*) \pi(s^*)$$
$$= \pi(s) \pi(t) a \varepsilon_{t^*} \pi(t^* s^*)$$
$$= \pi(s) \pi(t) \varepsilon_{t^*} a \pi(t^* s^*)$$
$$= \pi(st) \varepsilon_{t^*} a \pi(t^* s^*)$$
$$= \alpha_{st}(a)$$

Therefore, $\alpha$ is a partial action of $S$ on $\mathsf{A}$.

**Remark 2** As in the previous case, it is easy to see *that $\alpha^\pi$ and $\alpha^{\pi'}$* are equivalent partial actions if *$\pi$ and $\pi'$* are equivalent representations.

Given an inverse semigroup $S$, let $\pi : S \to \mathsf{B}$ be a partial representation of $S$ into $\mathsf{B}$, and let $J$ be the ideal of $\mathsf{B}$ generated by all elements of the form $a\pi(s) - a\pi(t)$ where $a \in \mathsf{B}$ and $s, t \in S$ are such that $s \leq t$. Define the map $\widetilde{\pi} : S \to \dfrac{\mathsf{B}}{J}$ by $\widetilde{\pi}(s) = \Phi(\pi(s))$, where $\Phi : \mathsf{B} \to \dfrac{\mathsf{B}}{J}$ is the quotient map. Obviously, $\widetilde{\pi}$ is a partial representation of $S$ on $\dfrac{\mathsf{B}}{J}$.

**Proposition 2.5** Let $\pi : S \to \mathsf{B}$ be a partial representation of $S$ into $\mathsf{B}$ and suppose that $\widetilde{\pi} : S \to \dfrac{\mathsf{B}}{J}$ is the partial representation mentioned above. Consider the subalgebra $\mathsf{A} \subseteq \dfrac{\mathsf{B}}{J}$ and partial action $\alpha^{\widetilde{\pi}}$ as in Lemma 3.4. Then the map $\phi_{\widetilde{\pi}} : \mathsf{A} \times_{\alpha^{\widetilde{\pi}}} S \ni \sum_{s \in S} a_s \delta_s \mapsto \sum_{s \in S} a_s \widetilde{\pi}(s) \in \dfrac{\mathsf{B}}{J}$ is a homomorphism of $\mathsf{K}$-algebras such that $\phi_{\widetilde{\pi}} \circ \pi_{\alpha^{\widetilde{\pi}}} = \widetilde{\pi}$.

*Proof.* Since $a_s \delta_s$ ($s \in S$) span $L$, it is enough to show that

$$\phi_{\widetilde{\pi}}(a_s \delta_s \cdot b_t \delta_t) = \phi_{\widetilde{\pi}}(a_s \delta_s) \cdot \phi_{\widetilde{\pi}}(b_t \delta_t).$$

By the definition of $\phi_{\widetilde{\pi}}$, we have

$$\phi_{\tilde{\pi}}(a_s\delta_s \cdot b_t\delta_t) = \phi_{\tilde{\pi}}(\alpha_s^{\tilde{\pi}}(\alpha_{s^*}^{\tilde{\pi}}(a_s)b_t)\delta_{st}) = \alpha_s^{\tilde{\pi}}(\alpha_{s^*}^{\tilde{\pi}}(a_s)b_t)\tilde{\pi}(st).$$

Since $\alpha_s^{\tilde{\pi}}(\alpha_{s^*}^{\tilde{\pi}}(a_s)b_t) \in X_s$ and $\varepsilon_s$ is the identity element of $X_s$,

$$\begin{aligned}
\alpha_s^{\tilde{\pi}}(\alpha_{s^*}^{\tilde{\pi}}(a_s)b_t)\tilde{\pi}(st) &= \alpha_s^{\tilde{\pi}}(\alpha_{s^*}^{\tilde{\pi}}(a_s)b_t)\varepsilon_s\tilde{\pi}(st) \\
&= \alpha_s^{\tilde{\pi}}(\alpha_{s^*}^{\tilde{\pi}}(a_s)b_t)\varepsilon_s\tilde{\pi}(s)\tilde{\pi}(t) \\
&= \tilde{\pi}(s)(\alpha_{s^*}^{\tilde{\pi}}(a_s)b_t)\tilde{\pi}(s^*)\tilde{\pi}(s)\tilde{\pi}(t) \\
&= \tilde{\pi}(s)(\alpha_{s^*}^{\tilde{\pi}}(a_s)b_t)\varepsilon_{s^*}\tilde{\pi}(t) \\
&= \tilde{\pi}(s)\alpha_{s^*}^{\tilde{\pi}}(a_s)b_t\tilde{\pi}(t) \\
&= \tilde{\pi}(s)\tilde{\pi}(s^*)a_s\tilde{\pi}(s)b_t\tilde{\pi}(t) \\
&= a_s\tilde{\pi}(s)b_t\tilde{\pi}(t) = \phi_{\tilde{\pi}}(a_s\delta_s)\phi_{\tilde{\pi}}(b_t\delta_t).
\end{aligned}$$

For an inverse semigroup $S$, let $\Pr(S)$ be the universal semigroup generated by the set of generators [s] chosen from a fixed set having as many elements as $S$ subject to the relations
  *(i)* $[s^*][s][t]=[s^*][st]$,
  *(ii)* $[s][t][t^*]=[st][t^*]$,
  *(ii)* $[s][s^*][s]=[s]$,
for all $s, t \in S$. We will refer the reader to [1] for more details about $\Pr(S)$.

**Definition 2.6** Let $S$ be an inverse semigroup. We will denote the semigroup algebra $K\Pr(S)$ by $K_{par}(S)$, and we will call it partial semigroup algebra.

Obviously, the canonical map $\iota_S : S \ni s \mapsto [s] \in K_{par}(S)$ is a partial homomorphism of $S$.

**Remark 3** Suppose that $\pi : S \to \mathsf{B}$ is a partial homomorphism of $S$ in $\mathsf{B}$. Then by [1, Proposition 2.20], there exists a unique semigroup homomorphism $\pi^* : \Pr(S) \to \mathsf{B}$ such that $\pi^* \circ \iota_S = \pi$. Now define $\psi : K_{par}(S) \ni \Sigma_{s \in \Pr(S)} a_s s \mapsto \Sigma_{s \in \Pr(S)} a_s \pi^*(s) \in \mathsf{B}$. Obviously, $\psi$ is a homomorphism and $\psi \circ \iota_S = \pi$. On the other hand, if $\psi : K_{par}(S) \to \mathsf{B}$ is an algebra homomorphism, then $\psi \circ \iota_S$ is a partial representation of $S$ in $\mathsf{B}$.

**Theorem 2.7** Let $\iota_S$ be the partial representation as mentioned above. Considering the partial representation $\tilde{\iota}_S$ and the subalgebra $\mathsf{A}$ of $\dfrac{K_{par}(S)}{J}$ generated by all elements $\tilde{\varepsilon}_s = \tilde{\iota}_S(s)\tilde{\iota}_S(s^*)$ defined in Proposition 2.5, we will prove that $\phi_{\tilde{\iota}_S}$ is an isomorphism.

Proof. Hereafter for simplicity, we will denote $\tilde{\iota}_S$ again by $\iota_S$. By Lemma 2.1, $\pi_{\alpha^{\iota_S}} : S \ni s \mapsto \tilde{\varepsilon}_s\delta_s \in \mathsf{A} \times_{\alpha^{\iota_S}} S$ is a partial representation of $S$ on $\mathsf{A} \times_{\alpha^{\iota_S}} S$. Now, by Remark 3, there exists $\psi : K_{par}(S) \to \mathsf{A} \times_{\alpha^{\iota_S}} S$ such that $\psi \circ \iota_S = \pi_{\alpha^{\iota_S}}$, i.e $\psi([s]) = \tilde{\varepsilon}_s\delta_s$. For $s \in S$, we have

$$\phi_{\iota_S} \circ \psi([s]) = \phi_{\iota_S}(\widetilde{\varepsilon}_s \delta_s) = \widetilde{\varepsilon}_s \iota_S(s) = [s][s^*][s] = [s].$$

Consequently, $\phi_{\iota_S} \circ \psi = id_{K_{par}(S)}$, so, $\phi_{\iota_S}$ is a surjective homomorphism. On the other hand, by the definition $X_s$, for each $s \in S$, is spanned by the element of the form $\widetilde{\varepsilon}_s \widetilde{\varepsilon}_{t_1} \ldots \widetilde{\varepsilon}_{t_n}$. For $s \in S$, we have

$$\widetilde{\varepsilon}_s \delta_s \widetilde{\varepsilon}_{s^*} \delta_{s^*} = \alpha_s^{\iota_S}(\alpha_{s^*}^{\iota_S}(\widetilde{\varepsilon}_s)\widetilde{\varepsilon}_{s^*})\delta_{ss^*}$$
$$= \alpha_s^{\iota_S}(\iota_S(s^*)\widetilde{\varepsilon}_s \iota_S(s)\widetilde{\varepsilon}_{s^*})\delta_{ss^*}$$
$$= \alpha_s^{\iota_S}([s^*][s][s^*][s]\widetilde{\varepsilon}_{s^*})$$
$$= \alpha_s^{\iota_S}(\widetilde{\varepsilon}_{s^*})\delta_{ss^*}$$
$$= [s]\widetilde{\varepsilon}_{s^*}[s^*] = \widetilde{\varepsilon}_s \delta_{ss^*} = \widetilde{\varepsilon}_s \delta_1.$$

Thus,

$$\psi \circ \phi_{\iota_S}((\widetilde{\varepsilon}_s \widetilde{\varepsilon}_{t_1} \ldots \widetilde{\varepsilon}_{t_n})\delta_s) = \psi(\widetilde{\varepsilon}_s \widetilde{\varepsilon}_{t_1} \ldots \widetilde{\varepsilon}_{t_n}[s])$$
$$= \psi([s][s^*][t_1][t_1^*]\ldots[t_n][t_n^*][s])$$
$$= \widetilde{\varepsilon}_s \delta_s \widetilde{\varepsilon}_{s^*} \delta_{s^*} \ldots \widetilde{\varepsilon}_{t_n} \delta_{t_n} \widetilde{\varepsilon}_{t_n} \delta_{t_n^*} \widetilde{\varepsilon}_s \delta_s$$
$$= \widetilde{\varepsilon}_s \delta_1 \widetilde{\varepsilon}_{t_1} \delta_1 \ldots \widetilde{\varepsilon}_{t_n} \delta_1 \widetilde{\varepsilon}_s \delta_s$$
$$= \widetilde{\varepsilon}_s \widetilde{\varepsilon}_{t_1} \ldots \widetilde{\varepsilon}_{t_n} \delta_s.$$

Thus, $\psi \circ \phi_{\iota_S}$ is also the identity map, by consequence, $\phi_{\iota_S}$ is an isomorphism.